\documentstyle[12pt]{article}
\makeatletter
\begin{document}
\thispagestyle{empty}
\pagestyle{empty}
\parskip = 10pt
\newcommand{\namelistlabel}[1]{\mbox{#1}\hfil}
\newenvironment{namelist}[1]{%
\begin{list}{}
{
\let\makelabel\namelistlabel
\settowidth{\labelwidth}{#1}
\setlength{\leftmargin}{1.1\labelwidth}
}
}{%
\end{list}}
\def\theequation{\thesection.\arabic{equation}}
\newtheorem{theorem}{\bf Theorem}[thetheorem]
\newtheorem{corollary}{\bf Corollary}[thecorollary]
\newtheorem{remark}{\bf Remark}[theremark]
\newtheorem{lemma}{\bf Lemma}[thelemma]
\def\thetheorem{\arabic{theorem}}
\def\theremark{\thesection.\arabic{remark}}
\def\thecorollary{\thesection.\arabic{corollary}}
\def\thelemma{\thesection.\arabic{lemma}}
\newcommand{\nc}{\newcommand}
\newcommand{\bsp}{\begin{sloppypar}}
\newcommand{\esp}{\end{sloppypar}}
\newcommand{\be}{\begin{equation}}
\newcommand{\ee}{\end{equation}}
\newcommand{\beanno}{\begin{eqnarray*}}
\newcommand{\inp}[2]{\left( {#1} ,\,{#2} \right)}
\newcommand{\dip}[2]{\left< {#1} ,\,{#2} \right>}
\newcommand{\disn}[1]{\|{#1}\|_h}
\newcommand{\pax}[1]{\frac{\partial{#1}}{\partial x}}
\newcommand{\tpar}[1]{\frac{\partial{#1}}{\partial t}}
\newcommand{\xpax}[2]{\frac{\partial^{#1}{#2}}{\partial x^{#1}}}
\newcommand{\pat}[2]{\frac{\partial^{#1}{#2}}{\partial t^{#1}}}
\newcommand{\ntpa}[2]{{\|\frac{\partial{#1}}{\partial t}\|}_{#2}}
\newcommand{\xpat}[2]{\frac{\partial^{#1}{#2}}{\partial t \partial x}}
\newcommand{\npat}[3]{{\|\frac{\partial^{#1}{#2}}{\partial t^{#1}}\|}_{#3}}
\newcommand{\eeanno}{\end{eqnarray*}}
\newcommand{\bea}{\begin{eqnarray}}
\newcommand{\eea}{\end{eqnarray}}
\newcommand{\ba}{\begin{array}}
\newcommand{\ea}{\end{array}}
\newcommand{\nno}{\nonumber}
\newcommand{\dou}{\partial}
\newcommand{\bc}{\begin{center}}
\newcommand{\ec}{\end{center}}
\newcommand{\bb}{\mbox{\hspace{.25cm}}}
\nc{\benu}{\begin{enumerate}}
\nc{\eenu}{\end{enumerate}}
\nc{\bth}{\begin{theorem}}
\nc{\eth}{\end{theorem}}
\nc{\bpr}{\begin{prop}}
\nc{\epr}{\end{prop}}
\nc{\blem}{\begin{lemma}}
\nc{\elem}{\end{lemma}}
\nc{\bcor}{\begin{corollary}}
\nc{\ecor}{\end{corollary}}
\newcommand{\R}{I\!\!\!R}
\nc{\Hy}{I\!\!\!H}
\nc{\C}{I\!\!\!C}
\newcommand{\la}{\lambda}
\newcommand{\s}{\sinh}
\newcommand{\co}{\cosh}
\newcommand{\vl}{V_{\la}}
\nc{\ga}{\gamma}
\nc{\vg}{V_{\ga}}
\nc{\T}{\Theta}
\nc{\f}{\frac}
\nc{\rw}{\rightarrow}
\nc{\om}{\omega}
\nc{\Om}{\Omega}
\nc{\al}{\alpha}
\nc{\qed}{\hfill \rule{2.5mm}{2.5mm}}
\setcounter{section}{0}
\title{A Pinching constant for Harmonic Manifolds}                 
\author{ K.Ramachandran and A.Ranjan}
\date{  }
\maketitle
\begin{abstract}
In this note we shall show that the sectional curvature of a
harmonic manifold is bounded on both sides. In fact we shall 
give a pinching constant for all harmonic manifolds. We shall use the
imbedding theorem for harmonic manifolds proved by Z.I.Szabo 
and the description of screw lines in hilbert spaces to prove the result.
\end{abstract}
\section{Introduction}
\label{sec-intro}
\setcounter{equation}{0}
Let $(M,g)$ be a riemannian manifold. It is well known that
radial harmonic functions, i.e solutions to $\Delta f=0$ which
depends only on the geodesic distance $r(x,.)$ , exists on $M$
only when the density funtion $\omega_{p}=\sqrt{| \mbox{det}(g_{ij})|}$
in a normal coordinate neibourhood around each point $p$ depends only on
the geodesic distance $r(p,.)$. 
A riemannian manifold is said to be harmonic if the density funtion 
satisfies the above radial property. The only known examples of Harmonic
spaces (for quite some time ) were the rank one symmetric spaces. The 
Lichnerowicz conjecture asserted that these are the only ones. In
\cite{Sza} Szabo proves that any compact harmonic manifold with finite
fundamental group is rank one symmetric. 
Later Damek and Ricci \cite{DR} gave examples of noncompact harmonic 
manifolds which are not rank one symmetric. Apart from these there are
no known examples of harmonic spaces. The sectional curvature of a harmonic
space was known to be bounded \cite{Bes}, our result gives the pinching
constant.  
\newpage
\section{The main Result}
\label{sec-result}
\setcounter{equation}{0}
Let $(M,g)$ be a riemannian manifold. Let $(x_1,...,x_n)$ be a normal
co-ordinate neibourhood around a point $p \in M$. The function 
$$\omega_p = \sqrt{| \mbox{det}(g_{ij}|}$$
is the volume density of $(M,g)$. The density function in polar co-ordinates
$(r_p,\phi)$ is then given by $\theta_p = r^{n-1}\omega_p$, where $r_p$ is  
the geodesic distance from $p$ and $\phi$ is a point on the unit sphere.
A riemannian manifold is said to be harmonic if the density function
$\theta_p$ is a funtion of the geodesic distance $r(p,.)$ alone. 
Besse \cite{Bes} constructed isometric imbeddings of compact harmonic
manifolds into their eigen spaces. In \cite{Sza} Szabo generalised Besse's
imbedding theorem. We shall state the generalised version.
Consider a $C^1$ function $h : \R_+ \rightarrow \R$ with $h'(0) = 0$ such that
$h,h' \in L^{2}_\theta(\R)$. We define the map 
$$\Phi_h : M \rightarrow L^{2}(M) \mbox{ \ by }$$
$$\Phi_h(p) = h_p,\;\; where\;\;  h_p(y) = h(r(p,y)).$$
\begin{theorem}
Let $(M,g)$ be a harmonic manifold. Then

\noindent
1. For any function $h$ as above the map 
$$ \Phi : M \rightarrow L^2(M) ;\;\; p \mapsto h_p$$
where $h_p(y) = h(r(p,y))$ is an isometric immersion of the harmonic space $M$
into a sphere of $L^2(M)$.

\noindent
2. The geodesics of $\Phi(M)$ are congruent screw lines in the space $L^2(M)$.
By screw lines in $L^2(M)$ we mean a rectifiable continuous curve $r(s)$
parametrized by arclength $s$ for which the distance $d( r(t) , r(s) )$ in
$L^2(M)$ depends only on the arclength $t-s$ for any two points $r(t) ,
r(s)$. \hfill{$\Box$} 
\end{theorem}

\hspace*{0.5cm}
we shall first study the second fundamental form of a harmonic manifold.
Neumann and Shoenberg \cite{NS} constructed for any screw line $r(t)$ in a
hilbert space $H$ a continuous one parameter family of unitary operators
$U(t) = e^{tX},\ X$ a skew symmetric operator such that $r(t) = e^{tX}v, v
= r(0)$. hence the general equation of a screw line is:
$$
\gamma(t)=e^{tX}v \;\mbox{for a skew symmetric operator}\; X.
$$
Let $B$ be the second fundamental form of $M$. As
$\gamma$ is a geodesic,
\begin{eqnarray*}
B(\gamma^{\prime},\gamma^{\prime})&=& \frac{d^2 \gamma}{d t^2} 
= X^2 e^{tX} v\\
\mbox{hence} \; |B(\gamma^{\prime},\gamma^{\prime})|^2 &=&
|X^2\;e^{tX}v|^2, \mbox{which is a constant}.
\end{eqnarray*}
If $\sigma$ is another geodesic congruent to $\gamma$, then
$\gamma=U\sigma$, where $U$ is a unitary operator. Hence
$B(\gamma',\gamma') = UB(\sigma',\sigma')$, hence
$|B(\gamma',\gamma')|^2 = |B(\sigma',\sigma')|^2$. The next
result is a general result on bilinear forms.
\begin{lemma}
Let $B : V\times V \rightarrow W$ be a symmetric bilinear form,
where $V,W$ are innerproduct spaces over $I\!\!R$,such that
$\|B(u,u)\|= c\|u\|^2$ ,where $c$ is a constant, then
$$ 2\|B(u,v)\|^2 = c^2( \|u\|^2+ \|v\|^2+ 2\inp{u}{v}^2 )-
\inp{B(u,u)}{B(v,v)}$$ .
\end{lemma}
{\bf Proof:}
$4B(u,tv)= B(u+tv,u+tv)-B(u-tv,u-tv)$, taking norms one gets
\begin{eqnarray*}
16t^2 \|B(u,v)\|^2 &=& \|B(u+tv,u+tv)\|^2 + \|B(u-tv,u-tv)\|^2\\
&-&2 \inp{B(u+tv,u+tv)}{B(u-tv,u-tv)}.
\end{eqnarray*}
Comparing coefficients of $t^2$ gives the result.\\
\begin{theorem}     
For any harmonic manifold $(M,g)$ the sectional curvature $K_M$ satisfies
$$ -2c^2 \leq K_M \leq c^2,\;\;\mbox{where c is a suitable constant.}$$
\end{theorem}
{\bf Proof:}
Embed the harmonic manifold $(M,g)$ into $L^2(M)$ such that all the
geodesics are congruent screw lines in $L^2(M)$. Let $B$ be the second
fundamental form of $M$, then $\|B(u,u)\|$ = constant($c$, say), for all
unit vectors $u$ at all points. Let $u,v$ be orthonormal unit tangent
vectors to $M$ at a point. Let P be the plane generated by $u,v$. The
sectional curvature is then,
$$ K(P) = \inp{B(u,u)}{B(v,v)}-\|B(u,v)\|^2$$.
using the above lemma one gets
$$K(P)=\frac{3}{2} \inp{B(u,u)}{B(v,v)}-\frac{1}{2}c^2$$
Schwartz inequality gives $ -c^2 \leq \inp{B(u,u)}{B(v,v)} \leq c^2$ which
gives 
$$ -2c^2 \leq K(P) \leq c^2.$$ 
The inequality on the left holds iff $B(u,u)=-B(v,v)$ and the one on the 
right holds iff $B(u,u)=B(v,v).$
\bcor
Let $M \subseteq S^n$ have all geodesics congruent screw lines in
$I\!\!R^{n+1}$. Suppose $codim M (\ in\ \R^{n+1}\ )=1,2$ then $M$ must be
of constant curvature. 
\ecor 
{\bf Proof:} 
Let $codim M =1$. In this case $dimW =1$ in the above lemma.
Let $\underline{e}$ be a unit vector in $W$, then we have $B(u,v)=\pm
\inp{u}{v} \underline{e}$. Hence $K(P)=c^2$ a constant. 
If $codimM =2$ we argue as follows. Let $p$ be any point of $M$. $p$ is normal
to $M$, let $\underline{n}$ be the other normal to $M$. Then for any
tangent vector $u$ to $M$ $\inp{B(u,u)}{p}=\inp{\gamma''(0)}{\gamma(0)}$
where $\gamma(t)$ is a geodesic at $p$. But $\inp{\gamma''(0)}{\gamma(0)}$
is a constant independent of the point $p$ since all geodesics are congruent
screw lines. Similarly $\inp{B(u,u)}{\underline{n}}$ is also a constant
independent of the point $p$. Again $M$ has constant curvature.
\bibliographystyle{plain}

\noindent
Department of Mathematics,\\    
I.I.T Powai, Mumbai 400 076,\\ 
India.\\
email: kram@ganit.math.iitb.ernet.in\\
\hspace*{1.25cm}aranjan@ganit.math.iitb.ernet.in\\
\end{document}